\documentclass[12pt]{amsart} 
\usepackage{amsthm, amssymb, amsmath,verbatim,float,graphicx} 
\usepackage[dvipsnames]{xcolor}
\usepackage{fullpage} 
\usepackage[latin1]{inputenc}
\numberwithin{equation}{section}
\newcommand{\lra}{\longrightarrow} 
\newcommand{\ra}{\rightarrow} 
\newcommand{\conic}{\mathcal K} 
\newcommand{\Rico}{\mathcal R} 
\newcommand{\Y}{\mathcal Y} 
\newcommand{\ry}{\mathrm y} 
\newcommand{\rz}{\mathrm z} 
\newcommand{\pline}{\mathbb L} 
\newcommand{\Sym}{\text{Sym}}
\newcommand{\SG}{\mathfrak G}
\newcommand{\U}{\mathbf U}
\newcommand{\ff}{\varphi} 
\newcommand{\He}{\text{He}}
\newcommand{\FF}{\mathcal F} 
\renewcommand{\P}{\mathbb P} 
\newcommand{\ux}{\mathbf x} 
\newcommand{\bA}{\mathbb A} 
\newcommand{\bB}{\mathbb B} 
\newcommand{\bC}{\mathbb C} 
\newcommand{\bD}{\mathbb D} 
\newcommand{\bE}{\mathbb E} 
\newcommand{\bF}{\mathbb F} 
\newcommand{\complex}{\mathbb C} 

\newcommand{\ltr}{\text{\sc ltr}}

\newcommand{\hltr}{\text{\sc h-ltr}}

\newcommand{\pasc}[6]{\left\{ \begin{array}{ccc} #1 & #2 & #3\\ #4 &
     #5 & #6 \end{array} \right\}}
\newtheorem{Theorem}{Theorem}[section]
\newtheorem{Lemma}[Theorem]{Lemma}
\newtheorem{Proposition}[Theorem]{Proposition}

\begin{document} 

\title{On the geometry of the ricochet configuration} 
\author{Jaydeep Chipalkatti} 
\maketitle 

\bigskip 

\parbox{17cm}{ \small
{\sc Abstract:} This paper is a study of the so-called `ricochet
configuration' (or $R$-configuration) which arises in the context of
Pascal's theorem. We give a geometric proof of the fact that a specific pair of Pascal lines
is coincident for a sextuple in $R$-configuration.  We calculate the
symmetry group of a generic $R$-configuration, as well as the degree of the subvariety 
$\Rico \subseteq \P^6$ of all such configurations. 
We also determine the $SL(2)$-equivariant defining equations for 
$\Rico$, and show that it is an ideal-theoretic complete
intersection of two invariant hypersurfaces.}  

\bigskip


\bigskip 

\tableofcontents

\section{Introduction} 
The `ricochet configuration' is a specific arrangement of six points
on a conic which arises in the context of Pascal's theorem. It was
discovered by the author in ~\cite{JC1}. 
We recall some of the background below for ease of reading. 

\subsection{} 
Let $\conic$ denote a nonsingular conic in the complex projective
plane. Consider six distinct points $A,B,C,D,E,F$ on $\conic$,
arranged into an array 
$\left[ \begin{array}{ccc} A & B & C \\ F & E & D \end{array}
\right]$. Then Pascal's theorem says that the three cross-hair intersection points 
\[ AE \cap BF, \quad BD \cap CE, \quad AD \cap CF\] 
(corresponding to the three minors of the array) are collinear. 

\begin{figure}[H] 
\includegraphics[width=8cm]{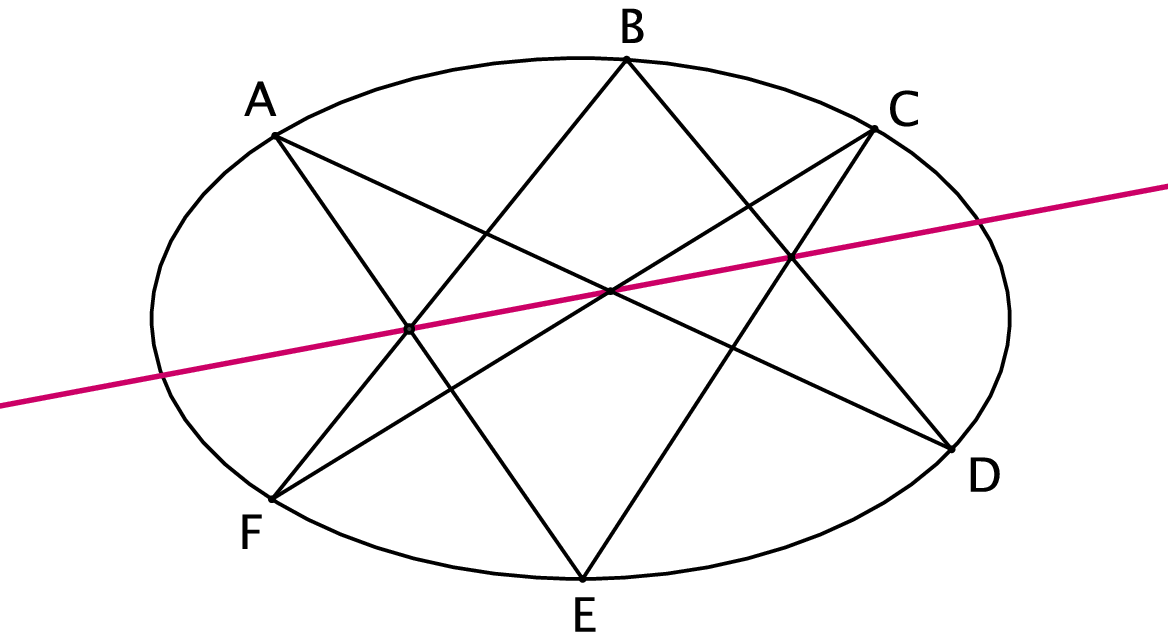}
\caption{Pascal's theorem} 
\end{figure}  

The line containing them is called the Pascal line, or just the Pascal, of
the array;  we will denote it by $\pasc{A}{B}{C}{F}{E}{D}$. It is easy to see that the Pascal remains unchanged if we permute the rows or the
columns of the array; for instance 
\[ 
\pasc{A}{B}{C}{F}{E}{D}, \quad \pasc{F}{E}{D}{A}{B}{C}, \quad 
\pasc{E}{D}{F}{B}{C}{A} \] 
all denote the same line. 

Any essentially different arrangement of the same sextuple of points, say 
$\pasc{E}{A}{C}{B}{F}{D}$, corresponds \emph{a priori} to a different
line. Hence we have a total of $\frac{6!}{2!3!} = 60$ notionally distinct Pascals. It is a
theorem due to Pedoe~\cite{Pedoe1}, that these $60$ lines are pairwise
distinct for a \emph{general} choice of the initial sextuple. In other
words, there must be something geometrically special about the
sextuple if some of its Pascals are to coincide. 

The main theorem on~\cite[p.~12]{JC1} characterises all such
special situations. It says that if some of the Pascals coincide,
then the sextuple must either be
in~\emph{involutive configuration}, or in~\emph{ricochet configuration}. We will describe both
of these below. The first is very classical
(cf.~\cite[\S 260]{Salmon}); whereas the second is probably not. To
the best of my knowledge, it had not previously appeared in the literature 
before it was discovered in the process of proving the theorem. 

\subsection{The involutive configuration} 
The sextuple $\Gamma = \{A, \dots, F \}$ is said be in involutive configuration, 
(or in involution for short), if there exists a point $Q$ in the plane with three
lines $L, L', L''$ through it such that 
\[ \Gamma = (L \cup L' \cup L'') \cap \conic. \] 
\begin{figure}[H] 
\includegraphics[width=12cm]{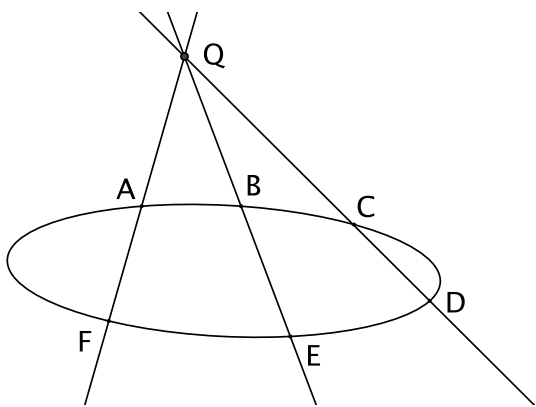}
\caption{The involutive configuration} 
\label{fig:inv} \end{figure}  
With points labelled as in the diagram, it turns out that the Pascals 
\[ \pasc{A}{B}{C}{F}{E}{D}, \quad 
\pasc{\textcolor{red}{F}}{B}{C}{\textcolor{red}{A}}{E}{D}, \quad 
\pasc{A}{\textcolor{red}{E}}{C}{F}{\textcolor{red}{B}}{D}, \quad 
\pasc{A}{B}{\textcolor{red}{D}}{F}{E}{\textcolor{red}{C}}, \] 
all coincide (see~\cite[p.~9]{JC1}). The pattern is straightforward: fix any column in the
first array and switch its entries to get another array. 
The common Pascal is the polar of $Q$ with
respect to the conic. 

\subsection{The ricochet configuration} 
\begin{figure}
\includegraphics[width=12cm]{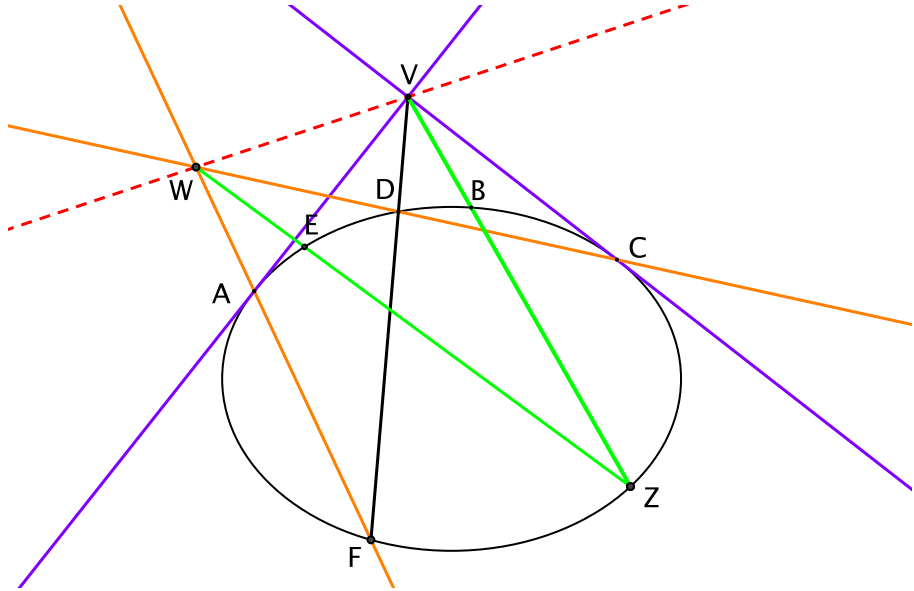}
\caption{The ricochet configuration} 
\label{fig:ricochet} 
\end{figure}

The construction in this case is rather more elaborate. 
Start with arbitrary distinct points $A,B,C,D$ on the conic. We
will define two more points $E$ and $F$ to complete the sextuple
(see Diagram \ref{fig:ricochet}). 

\begin{itemize} 
\item Draw tangents to the conic at $A$ and $C$. Let $V$ denote their 
  intersection point. 
\item Extend $VD$ so that it intersects the conic again at $F$. 
\item Let $W$ be the intersection point of $AF$ and $CD$. 
\item Now mark off $Z$ on the conic such that $V, B, Z$ are collinear,
and finally $E$ such that $W,Z,E$ are collinear. 
\end{itemize}

One can think of $B$ as a billiard ball which is struck by $V$ so that 
it bounces off the conic at $Z$, and gets redirected to $W$; hence the
name `ricochet'. 
For such a sextuple, the Pascals 
\begin{equation} 
\pasc{A}{\textcolor{OrangeRed}{B}}{C}{\textcolor{blue}{F}}
{\textcolor{OrangeRed}{E}}{\textcolor{blue}{D}}, \qquad 
\pasc{A}{\textcolor{OrangeRed}{E}}{C}{\textcolor{blue}{D}}
{\textcolor{OrangeRed}{B}}{\textcolor{blue}{F}} 
\label{coinc.pascals.rico} \end{equation} 
coincide (see~\cite[p.~10]{JC1}). The common Pascal is the line $VW$,
something which is not altogether obvious from the
diagram. It is \emph{prima facie} a little odd that the Pascal only
depends on $A, C, D$ and not on $B$. All of this will be clarified in section~\ref{section.doubleinv}. 

As mentioned above, the main result of~\cite{JC1} can be paraphrased
as saying that any sextuple for which some of the Pascals coincide
must fit into either Diagram~\ref{fig:inv} or 
Diagram~\ref{fig:ricochet}, up to a relabelling of points. 

\subsection{A summary of results} This paper is a study of the
algebro-geometric properties of the 
ricochet configuration (henceforth called the $R$-configuration). 
\begin{enumerate} 
\item 
The fact that the two Pascals in (\ref{coinc.pascals.rico}) coincide 
was proved by an inelegant brute-force calculation
in~\cite{JC1}. We will give a geometric proof in section~\ref{section.doubleinv}. 
\item 
In section~\ref{section.shufflegroup}, we determine the group of symmetries
of a generic $R$-configuration; it turns out to be the $8$-element
dihedral group. If one thinks of a sextuple as an element of 
\[ \Sym^6 \, \conic \simeq \Sym^6 \, \P^1 \simeq \P^6, \] 
then all sextuples in ricochet configuration form a $4$-dimensional subvariety $\Rico \subseteq
\P^6$. The symmetry group will be used to prove that $\Rico$ has degree $60$. 
\item 
The special linear group $SL(2,\complex)$ acts on the
projective plane by linear automorphisms in such a way that $\conic$
is stabilized (more on this in section~\ref{section.SLrepresentations} below). Since the $R$-configuration is
constructed synthetically, the subvariety $\Rico$ is 
stabilized by the induced action on $\Sym^6 \, \conic$. It 
follows that $\Rico$ must be defined by $SL(2)$-invariant homogeneous
equations; or in classical language, by the vanishing of certain
covariants of binary sextic forms. We will find such equations explicitly
in section~\ref{section.equiv.equations}. It turns out that $\Rico$ is defined by the
vanishing of two invariants, one each in degrees $6$ and $10$. 
\end{enumerate}

All the necessary background in projective geometry may be 
found in~\cite{Coxeter, Seidenberg, SempleKneebone}. We will
use~\cite{Smith_etal} as the standard reference for algebraic
geometry, but nothing beyond the most basic notions will be needed. 

\section{Preliminaries} 
Our entire set-up agrees with the one used in ~\cite[Ch.~3]{JC1}. We 
will recall only some of it below, and refer the reader to the earlier
paper for details. Section~\ref{section.doubleinv} is in any event entirely
geometric, and apart from section~\ref{section.involutions} on
involutions it does not need any of the algebraic preliminaries given
here. 

\subsection{} 
\label{section.SLrepresentations} 
For $m \ge 0$, let $S_m$ denote the vector space of homogeneous 
polynomials of degree $m$ in the variables $\ux = \{x_1, x_2\}$. In
classical language, elements of $S_m$ are the binary $m$-ics. Given $A
\in S_m$ and $B \in S_n$, their $r$-th transvectant will be denoted by 
$(A,B)_r$. It is a binary form of degree $m+n-2r$. 

We will use $\P^2 = \P S_2$ as our working projective plane; thus a nonzero
quadratic form $Q = a_0 \, x_1^2 + a_1 \, x_1 \, x_2 + a_2 \, x_2^2$
represents a point $[Q] \in \P^2$. Consider the Veronese imbedding 
\[ \P S_1 \stackrel{\phi}{\lra} \P S_2, \quad [u] \lra [u^2]. \] 
The image of $\phi$ will be our conic $\conic$. The point $[Q]$
lies on $\conic$, iff $Q$ is the square of a linear form. Thus 
$\conic$ is defined by the equation $a_1^2 = 4 \, a_0 \,
a_2$. Henceforth we will write $Q$ for $[Q]$ etc., if no confusion is
likely. We will sometimes use affine coordinates on $\conic \simeq \complex  
\cup \{\infty\}$, so that $\alpha  
\in \complex $ corresponds to $\phi(x_1 - \alpha \, x_2)$, and  
$\infty$ to $\phi(x_2)$.

The advantage of such a set-up is that the action of the special
linear group is naturally built into it. A matrix $M
= \left[ \begin{array}{cc} \alpha & \gamma \\ \beta &
                                                     \delta \end{array}
                                                 \right] \in SL(2, \complex)$ gives
an automorphism of $S_m$ defined by a linear change of variables 
$f(x_1, x_2) \ra f(\alpha \, x_1 + \beta \, x_2, \gamma \, x_1 +
\delta \, x_2)$; this in turn induces an automorphism of the projective space $\P S_m \simeq
\P^m$. The operation of transvection commutes with a linear change of
variables; in particular all the notions involving points and lines in $\P^2$, 
as well as polarities with respect to $\conic$ are expressible in
the language of transvectants. 

\subsection{Involutions} 
\label{section.involutions}
\begin{figure}
\includegraphics[width=6cm]{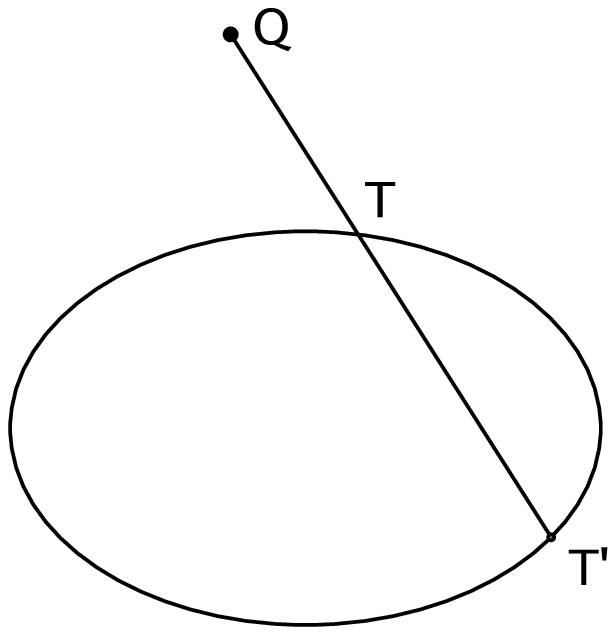}
\caption{$T \lra \sigma_Q(T)=T'$} 
\label{fig:inv1} 
\end{figure}

Every point $Q \in \P^2 \setminus \conic$ defines an involution (i.e.,
a degree $2$ automorphism) $\sigma_Q$ on $\conic$.  It takes a point $T$ to the
other intersection of $QT$ with $\conic$ (see
Diagram~\ref{fig:inv1}). In particular, $\sigma_Q(T)
= T$ exactly when $QT$ is tangent to $\conic$. 

\begin{figure}
\includegraphics[width=8cm]{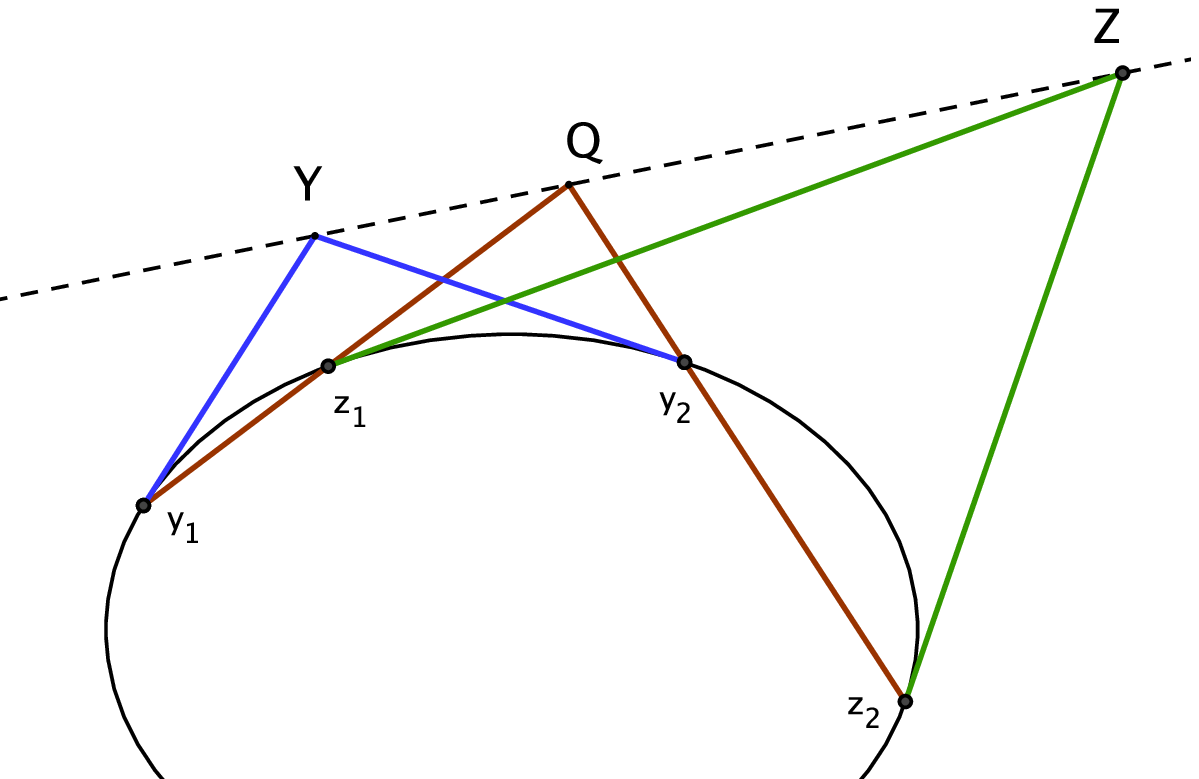}
\caption{$Y \lra \sigma_Q(Y)=Z$} 
\label{fig:inv2}
\end{figure}

Now let $Y$ be any point in $\P^2$, and let $\ry_1, \ry_2$ be the
intersection points of its polar with respect to $\conic$. (That is to say,
$\ry_i Y$ are tangent to the conic.) We define
$Z = \sigma_Q(Y)$ to be the pole of the line joining $\rz_1 = \sigma_Q(\ry_1),
\rz_2 = \sigma_Q(\ry_2)$. Thus $\sigma_Q$ extends to an involution of the entire
plane (see~Diagram~\ref{fig:inv2}). The points $Y, Q$ and $\sigma_Q(Y)$ are collinear. This has the 
consequence that if $\ell$ is a line passing through $Q$, then
$\sigma_Q(\ell) = \ell$ as a set. 
\subsection{Algebraic form of the $R$-configuration} We will express
the notion of an $R$-configuration in the language of
section~\ref{section.SLrepresentations}. 
Consider the set of letters $\ltr = \{\bA, \bB, \bC, \bD, \bE, \bF
\}$. Define a hexad to be an injective map 
$\ltr \stackrel{h}{\lra} \conic$, and write 
\[ A = h(\bA), \quad B = h(\bB), \quad \dots \quad F = h(\bF) \] 
for the corresponding distinct points on the
conic. Then $\Gamma = \text{image}(h) = \{A, \dots, F\}$ is the
associated sextuple. 
A hexad $h$ will be called an \emph{alignment} if the two rows of the
table 
\begin{equation} 
\begin{array}{c|c|c|c|c|r} A & B & C &  D &  E & F \\ \hline 
0 &  t & \infty & 1 & \frac{t-1}{t+1} &  -1 \end{array} 
\label{sequence.rico} \end{equation} 
are projectively isomorphic for some complex number $t$. In other
words, there should exist an automorphism of $\conic$ (or what is the
same, a fractional linear transformation of $\P^1$) which
takes $A, B, \dots, F$ respectively to $0, t, \dots, -1$. Since the points are required to be
distinct, we must have $t \neq 0, 1, \sqrt{-1}$. For later reference,
let $\Sigma(t)$ denote the sextuple corresponding to the second row of (\ref{sequence.rico}).

We will say that a sextuple $\Gamma$
is in $R$-configuration (or, it is an $R$-sextuple), if it
admits at least one alignment. To see that this definition agrees with the
geometric construction, choose coordinates on $\conic$ such
that $A,C,D$ respectively correspond to $0, \infty, 1$. This can
always be done by the fundamental theorem of projective
geometry. Using binary forms, 
\[ A = x_1^2, \quad C = x_2^2, \quad D = (x_1-x_2)^2.   \] 
Following the geometric construction, we get $V = x_1 \, x_2, F = (x_1
+x_2)^2$, and hence\footnote{As in \cite{JC1}, we will use $\Box$ to
  indicate a nonzero multiplicative scalar whose precise value is irrelevant.}
\[ W = (x_1 (x_1+x_2), x_2 \, (x_1 - x_2))_1 =
\Box \, (x_1^2-2 \, x_1 \, x_2-x_2^2). \] 

Now let $B = (x_1-t \, x_2)^2$ for some $t$. Then $Z = \sigma_V(B) =
\Box \, (x_1 + t \, x_2)^2$, and finally 
\[ E = \sigma_W(Z) = \Box \, (x_1 - \frac{t-1}{t+1} \, x_2)^2. \] 
This agrees exactly with (\ref{sequence.rico}). The 
ricochet $B \leadsto Z \leadsto E$ corresponds to $t
\ra -t \ra \frac{t-1}{t+1}$. Notice that $A,C,D,F$ is a harmonic quadruple, i.e., the cross-ratio 
$\langle A, C, D, F \rangle =-1$. Thus one can 
think of the $R$-configuration as a `fixed' harmonic quadruple, 
joined by a moving pair of points $B$ and $E$. 
It will be convenient to introduce the partition 
\begin{equation} 
\ltr = \underbrace{\{\bA,\bC,\bD,\bF\}}_{\hltr} \, \cup \, \{\bB,
\bE\},  
\label{hltr.defn} \end{equation} 
where $\hltr$ is to thought of as the `harmonic' subset of letters.

The fractional linear transformation 
\begin{equation} \ff(t) = \frac{t-1}{t+1}, 
\label{ff.formula} \end{equation} 
will appear many times below. Its inverse is given by $\ff^{-1}(t) =
\frac{1+t}{1-t}$. 

\subsection{Example} The table 
\[ 
\begin{array}{r|r|c|r|r|r} - \frac{1}{3} & 1 & 2 &  
\frac{1}{4} & \frac{1}{18} & - \frac{3}{2} \\ \hline 
0 &  4 & \infty & 1 & \frac{3}{5} &  -1 \end{array} \] 
is so arranged that the second row is $\Sigma(4)$, and $s \ra \frac{3s+1}{2-s}$ 
transforms the first row into the second. Hence
the first row (and of course, also the second) is in $R$-configuration. 
\subsection{} 
We identify the projective space $\P^6$ with $\P S_6$. A nonzero
binary sextic form $F$ will factor as $\prod\limits_{i=1}^6 (\alpha_i
\, x_1 - \beta_i \, x_2)$, and as such corresponds to the 
sextuple of points $\{\beta_i/\alpha_i: 1 \leqslant i \leqslant 6\}$
on $\P^1 \simeq \conic$. The points are
distinct if $F$ has no repeated linear factors. Hence the set of
sextuples of distinct points on $\conic$ can be identified with the complement of the discriminant
hypersurface in $\P^6$. 

Let $\Rico \subseteq \P^6$ denote the Zariski closure of the set of
all $R$-configurations; in other words, it is the Zariski closure of the union
of $SL(2)$-orbits of the sextic forms 
\[ G_t = x_1 \, x_2 \, (x_1 - x_2) \, (x_1 + x_2) \, (x_1 - t \, x_2) 
\, (x_1 - \ff(t) \, x_2), \] 
over all complex numbers $t \neq 0, 1, \sqrt{-1}$.  Since an
$R$-configuration is built from an arbitrary choice of $A, B, C,
D$ on the conic, $\Rico$ is an irreducible
$4$-dimensional rational projective variety. 

\section{The double involutions} 
\label{section.doubleinv} 
Let $\Gamma = \{A, \dots, F\}$ be in $R$-configuration. Our object is to show that both Pascals
in~(\ref{coinc.pascals.rico}) are equal to
the line $VW$ in~Diagram~\ref{fig:ricochet}. We will approach
the issue obliquely. 

\subsection{} 
\begin{figure}
\includegraphics[width=10cm]{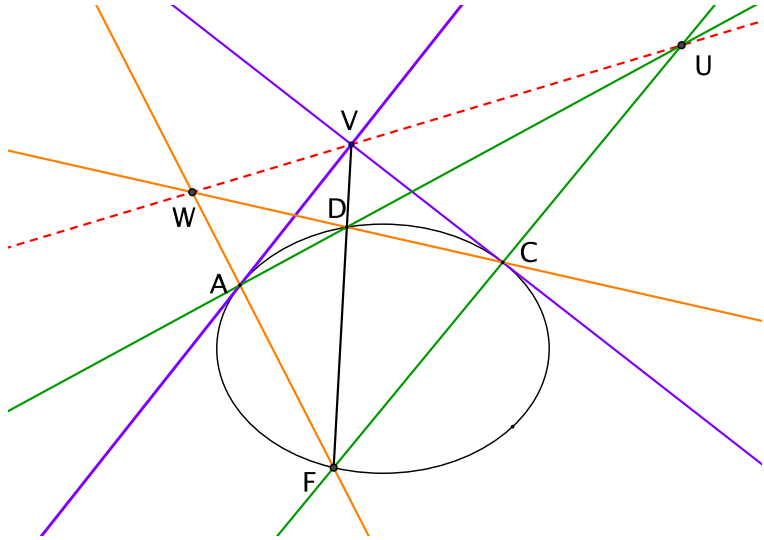}
\caption{The harmonic quadruple $\{A,C,D,F\}$} 
\label{fig:ricochet2}
\end{figure}

Diagram~\ref{fig:ricochet2} is a modified version of
Diagram~\ref{fig:ricochet}. The points
$A, C, D, F, V, W$ are exactly as before, but $B$ and $E$ are not yet
in the picture. The lines $AD, CF$ intersect in $U$. 

\begin{Lemma} \rm 
The points $U, V, W$ are collinear. 
\end{Lemma} 
\proof The involution $\sigma_V$ preserves the points $A, C$, and
interchanges $D, F$. Hence it takes $U = AD \cap CF$ to $W = AF \cap
CD$. Thus $U, V, W$ are collinear. \qed 

Let $\pline$ denote the line $UVW$. 

\begin{Lemma} \rm 
The automorphisms $\sigma_W \circ \sigma_V$ and $\sigma_V \circ
\sigma_U$ of $\conic$ are equal. 
\end{Lemma} 
\proof 
Since $\conic \simeq \P^1$, by the fundamental theorem of projective
geometry it will suffice to show that the two agree on three distinct
points. Now $\sigma_W \circ \sigma_V(A) = \sigma_W(A)= F$, and 
$\sigma_V \circ \sigma_U(A) = \sigma_V(D)= F$. Similarly, it is easy
to check that both maps send $C$ to $D$, and $D$ to $A$. \qed

\medskip 

Let $\psi: \conic \lra \conic$ denote this automorphism. 
For arbitrary points $B, E$ on the conic, consider the Pascal 
$\pasc{A}{B}{C}{F}{E}{D}$. By definition, it must pass
through the point $U = AD \cap CF$. As 
$B$ and $E$ move on the conic, the Pascal will pivot around $U$. We
should like to know under what conditions it will equal $\pline$. This is answered by the
next proposition. 

\begin{Proposition} \rm 
We have $\pasc{A}{B}{C}{F}{E}{D} = \pline$, exactly when $\psi(B) =
E$. 
\end{Proposition} 
\begin{figure}
\includegraphics[width=8cm]{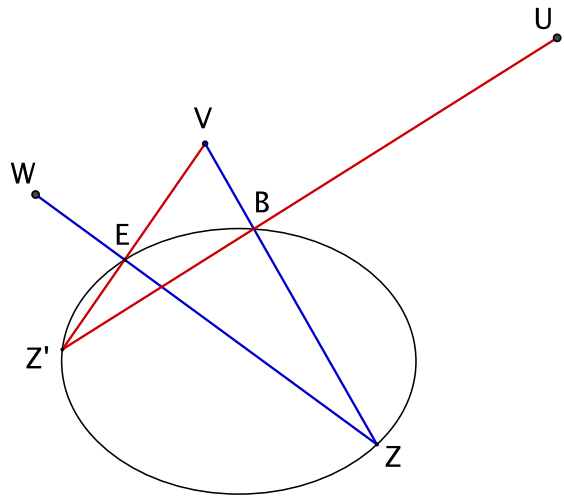}
\caption{$\sigma_W \circ \sigma_V(B)=E = \sigma_V \circ \sigma_U(B)$}
\label{fig:actionpsi} 
\end{figure}  
Diagram~\ref{fig:actionpsi} shows the action of $\psi = \sigma_W \circ
\sigma_V = \sigma_V \circ \sigma_U$. 
One can move from $B$ to $E$ either by a ricochet at $Z$ or at $Z'$. 
Let us assume the proposition for now, and deduce the equality of
Pascals. If $\psi(B) = E$, then $\psi(Z') = Z$. 
Applying the proposition with $Z'$ in place of $B$, we have 
\[ \pasc{A}{Z'}{C}{F}{Z}{D} = \pline. \] 
Now apply $\sigma_V$ to this equation. Since $\pline$ passes through
$V$, we have $\sigma_V(\pline) = \pline$ by
section~\ref{section.involutions}. But then 
\[ \pasc{\sigma_V(A)}{\sigma_V(Z')}{\sigma_V(C)}
{\sigma_V(F)}{\sigma_V(Z)}{\sigma_V(D)} = 
\pasc{A}{E}{C}{D}{B}{F} = \pline,  \] 
which is exactly what we wanted. 

\subsection{} 
It remains to prove the proposition. Given an arbitrary point $B$ on
the conic, we will define $\omega(B)$ such that 
\begin{equation} \pasc{A}{B}{C}{F}{\omega(B)}{D} = \pline. 
\label{eqn.omega.pascal} 
\end{equation} 
Afterwards we will prove that $\omega$ and $\psi$ are the same
morphism. One can define $\omega(B)$ in either of the following two ways; the
identity in~(\ref{eqn.omega.pascal}) is then simply the definition of the
Pascal. 
\begin{figure}[H] 
\includegraphics[width=8cm]{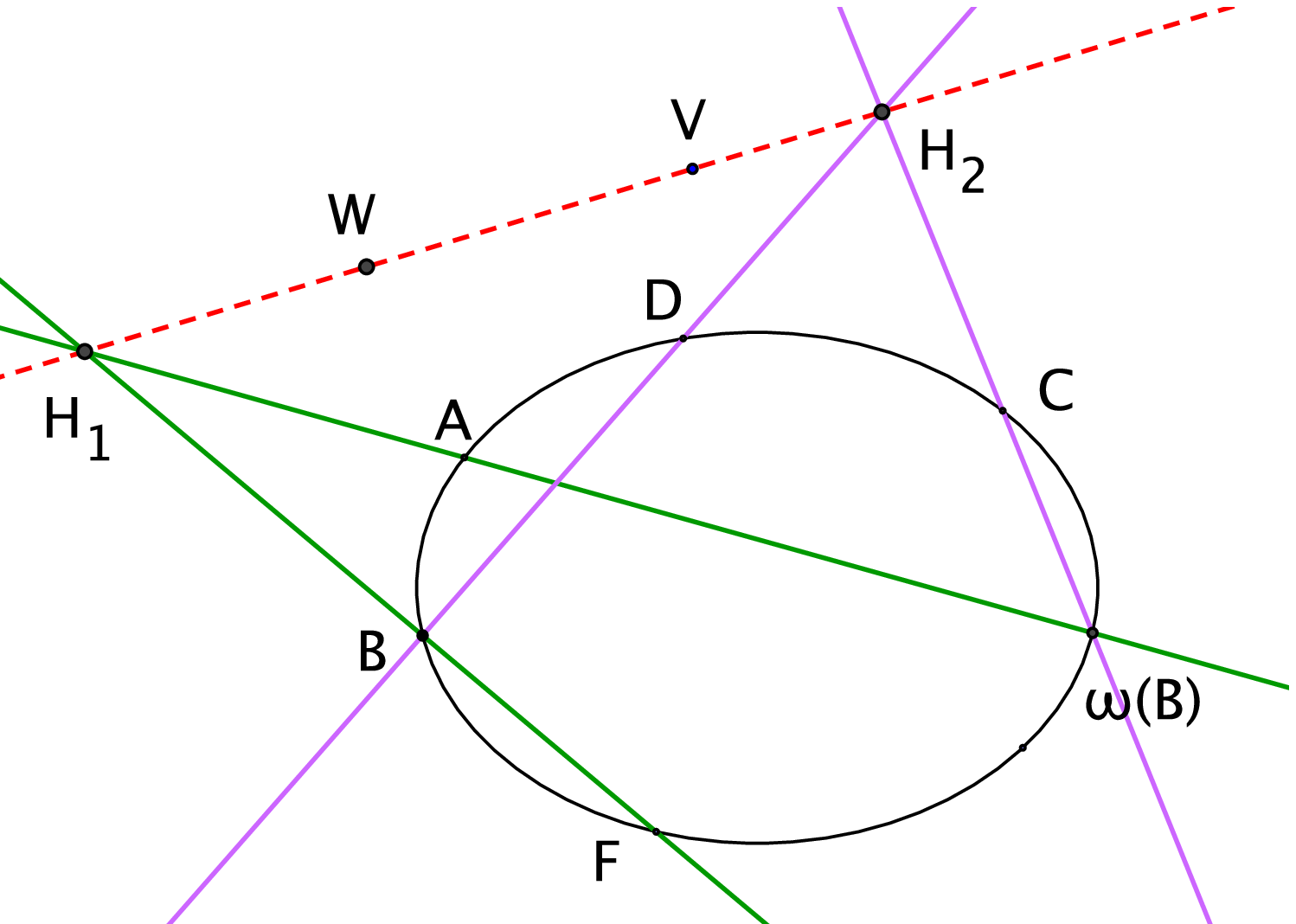}
\caption{$B \lra \omega(B)$}
\end{figure}  

\begin{enumerate} 
\item 
Intersect $BF$ with $\pline$ to get a point $H_1$, and define
$\omega(B)$ to be the other intersection of $A H_1$ with $\conic$. 
\item 
Intersect $BD$ with $\pline$ to get a point $H_2$, and define 
$\omega(B)$ to be the other intersection of $C H_2$ with $\conic$. 
\end{enumerate} 
This defines a morphism $\omega: \conic \lra \conic$, which is
bijective since the construction can be reversed to define
$\omega^{-1}$. One point should be clarified. Throughout this paper,
we have considered sextuples of distinct points only. However,
$\omega$ is defined for \emph{all} positions of $B$ on $\conic$, even
those which coincide with other points. For instance, if $B$ coincides with $D$,
then we interpret $BD$ as the tangent at $D$. 

Now observe that 
\begin{itemize} 
\item If $B = A$, then $H_1=W$ and $\omega(B) = F$. 
\item If $B = C$, then $H_2 = W$ and $\omega(B) = D$. 
\item 
If $B = D$, then $H_1 = V$ and $\omega(B) = A$ since 
$V$ lies on the tangent at $A$. 
\end{itemize} 
Thus $A, C, D$ are respectively mapped to $F, D, A$ by
$\psi$ as well as $\omega$, hence they must be the same morphism. This proves
the proposition. \qed 

The equality of $\psi$ and $\omega$ seems difficult to prove directly,
since their definitons are rather disparate. But the fundamental theorem of
projective geometry allows us to conclude the
argument by comparing their values only at three chosen points. In
summary, we have a geometric proof of the following theorem: 
\begin{Theorem} \rm 
If $\Gamma = \{A, \dots, F\}$ is a sextuple in $R$-configuration, then 
\[ \pasc{A}{B}{C}{F}{E}{D} = \pasc{A}{E}{C}{D}{B}{F}. \] 
\end{Theorem} 

In~\cite[p.~17]{JC1}, Gr{\"obner} basis computations are used to prove
that the converse of this theorem is also true, i.e., assuming that the
two Pascals coincide forces $\Gamma$ to be in $R$-configuration. It
would be interesting to have a geometric proof of this fact, but it is
not clear how to proceed. 

\section{The shuffle group and the degree of the ricochet locus} 
\label{section.shufflegroup} 
In this section we will determine the group of combinatorial
symmetries of a generic $R$-sextuple. This calculation will be of use 
in finding the degree of the variety $\Rico$. 
As in \cite{JC1}, let $\SG(X)$ denote the group of bijections $X \lra
X$ on a set $X$. 

\subsection{} 
Let $\Gamma = \Sigma(t)$ as in (\ref{sequence.rico}). Fix the 
alignment $h: \ltr \lra \Sigma(t)$ such that 
\[ \bA \ra 0, \quad \bB \ra t, \quad \bC \ra \infty, \quad \bD \ra 1,
\quad \bE \ra \frac{t-1}{t+1}, \quad \bF \ra -1. \] 

Consider the subgroup $H(t) \subseteq \SG(\ltr)$ consisting of elements $z$ such
that $h \circ z$ is also an alignment. In other words, $H(t)$
measures in how many ways the same sextuple can be seen to be in 
$R$-configuration. We may call it the shuffle group corresponding to
$t$. 

\begin{Lemma} \rm 
The elements 
\begin{equation} u = (\bA \, \bD \, \bC \, \bF), \quad 
v = (\bA \, \bD) (\bB \, \bE) (\bC \, \bF) 
\label{defns.uv} \end{equation} 
are in $H(t)$. 
\end{Lemma} 
By our convention, the $4$-cycle $u$ takes $\bA$ to $\bD$ etc. 
\proof 
The proof for $u$ is captured by the following table: 
\[ \begin{array}{c|c|r|c|r|r} 
1 & t & -1 & \infty & \frac{t-1}{t+1} & 0 \\ \hline 
0 & \frac{t-1}{t+1} & \infty & 1 & - \frac{1}{t} & -1 
\end{array} \] 
The hexad $h \circ u$ is given by $\bA \ra \bD \ra 1, \bB \ra \bB
\ra t$ etc, all of which is described by the first row. The fractional
linear transformation $s \ra \frac{s-1}{s+1}$ converts it into the
second row, which is $\Sigma(\frac{t-1}{t+1})$. Hence $h \circ u$ is
also an alignment, i.e., $u \in H(t)$. 

Similarly, the hexad $h \circ v$ is the first row of the table: 
\[ \begin{array}{c|c|r|c|r|r} 
1 & \frac{t-1}{t+1} & -1 & 0 & t & \infty \\ \hline 
0 & \frac{1}{t} & \infty & 1 & \frac{1-t}{1+t} & -1 
\end{array} \] 
The transformation $s \ra \frac{1-s}{1+s}$ converts it into the second
row, which is $\Sigma(\frac{1}{t})$. Hence $v \in H(t)$. \qed 

These two elements satisfy the relations $u^4 = v^2 = (u \, v)^2
=e$, hence the subgroup generated by them is the dihedral group with
$8$ elements. 

\subsection{} 
We already know that $\{A,C,D,F\}$ is a harmonic quadruple inside $\Sigma(t)$. But some 
other quadruple, say $\{B,D,C,E\}$, will be harmonic exactly when the
cross-ratio
\[ \langle B, D, C, E \rangle =\frac{2}{t^2+1}, \] 
is $-1, \frac{1}{2}$ or $2$. This can happen only for finitely many
values of $t$, and of course likewise for all such cases. 
Hence, $\{A,C,D,F\}$ is the only harmonic quadruple  
inside $\Sigma(t)$, for all but finitely many values of $t$ (that is
to say, for a `generic' $t$). 

\begin{Proposition} \rm 
For generic $t$, the group $H(t)$ is generated by $u$ and $v$. 
\end{Proposition} 
\proof By what has been said, every element in $H(t)$ must preserve
the subset $\hltr \subseteq \ltr$. This gives a morphism $f: H(t) \lra \SG(\hltr)$. 

Let $G \subseteq \SG(\hltr)$ denote the group of
permutations $\delta$ such that 
\[ \langle h \circ \delta(\bA), h \circ \delta(\bC), h \circ
\delta(\bD), h \circ \delta(\bF) 
\rangle =-1. \] 
Now $G$ is the $8$-element dihedral group generated by 
$u=(\bA \, \bD \, \bC \, \bF)$ and $v'=(\bA \, \bD) (\bC \,
\bF)$; this is a standard fact about the symmetries of the cross-ratio
and in particular those of a harmonic quadruple (see~\cite[Ch.~IV]{VeblenYoung}). 
We know,~\emph{a priori}, that the image of $f$ is contained in $G$. 
Now $f$ surjects onto $G$, since $f(v)=v'$. It is easy to check that $(\bB \, \bE)
\notin H(t)$, and hence $f$ is also injective. It follows that $f$ is
an isomorphism onto $G$, and thus $H(t)$ is generated by $u$ and $v$. \qed 

Let $H$ denote this group in the generic case. 

\subsection{} The group $H(t)$ may be
larger for special values of $t$. If $t = \sqrt{-3}$, then 
$\{B, D, C, E \}$ is also a harmonic quadruple, which allows more
possibilities for elements in $H(t)$. A routine computation shows that $H(\sqrt{-3})$ is
the $16$-element group generated by $u$ and $v$, together with the 
additional element $(\bA \, \bB) (\bC \, \bD) (\bE \, \bF)$. There are
several such special values of $t$, but we do not attempt to classify them. 

In general, an element of $H$ does not extend to an
automorphism of the entire conic. For instance, let $T$ denote the
intersection of the lines $AD, BE$. If $v$ were to extend to an
automorphism of $\conic$, it would have to coincide with the
involution $\sigma_T$, since both have identical actions on the four
points $A, B, D, E$. However this is a contradiction, since the
line $CF$ will not pass through $T$ for generic $t$. 

\subsection{} \label{examples.type} 
Now we will use the symmetry group to determine 
the degree of $\Rico$ as a projective subvariety in $\P^6$. 
If $z \in \conic$ is an arbitrary point, then $\{ \Gamma \in 
\P^6: z \in \Gamma \}$ is a hyperplane in $\P^6$. Since the 
degree of $\Rico$ is the number of points in its 
intersection with four general hyperplanes, we are reduced to the 
following question: Given a set of four general points $Z = \{z_1, \dots, z_4 
\} \subseteq \P^1 \simeq \conic$, find the number of $R$-sextuples 
$\Gamma$ which contain $Z$. 

Thus the degree of $\Rico$ can be understood in the following
intuitive way. Since the $R$-configuration has four degrees of
freedom, four general points on the conic 
will fit into only finitely many  $R$-configurations. 
We wish to know how many.\footnote{Since there are only finitely many values of $t$ for
which $H(t)$ is larger than $H$, the subclass of such
$R$-configurations has only three degrees of freedom. Hence a general
set of four points on the conic is not extendable to any such
configuration. This fact will play a role in the degree calculation
below.} The following two examples should capture the gist of the matter. 

Let $Z =\{2,3,5,7\}$, and assign them respectively to positions $\bA, \bC, \bD,
\bE$. This means that, in the table below 
\[ \begin{array}{r|r|c|r|c|r} 
A & B & C & D & E & F \\ \hline 
2 & \textcolor{red}{b} & 3 & 5 & 7 & \textcolor{red}{f} \\ 
0 & t & \infty & 1 & \ff(t) & -1 
\end{array} \] 
we want to find all pairs $(b,f)$ such that the second row is
projectively isomorphic to the third row for some $t$. The transformation $\mu(s)=
\frac{9s-4}{3s-2}$ takes $0,\infty,1$ respectively to $2,3,5$. Hence
$f=\mu(-1)=\frac{13}{5}$. Then $\varphi(t) = \mu^{-1}(7) =
\frac{5}{6}$, and hence $t = \varphi^{-1}(\frac{5}{6}) = 11$. Finally 
$b = \mu(11) = \frac{95}{31}$. Thus we have a unique pair $(b,f)$
which extends $Z$ to an $R$-sextuple. 

Now assign the same numbers to $\bA, \bB, \bD, \bE$, which leads to
the table: 
\[ \begin{array}{r|r|r|r|r|r} 
A & B & C & D & E & F \\ \hline 
2 & 3 & \textcolor{red}{c} & 5 & 7 & \textcolor{red}{f} \\ 
0 & t & \infty & 1 & \ff(t) & -1 
\end{array} \] 
As before, we are searching for all pairs $(c,f)$ such that the second and
the third rows are projectively isomorphic. 
Now $\nu(s) = \frac{(c-5) \, (s-2)}{3 \, (c-s)}$ takes $2, c, 5$
respectively to $0, \infty, 1$. Hence $t = \nu(3) = \frac{(c-5)}{3 \,
  (c-3)}$, which leads to $\ff(t) = \frac{t-1}{t+1} = \frac{c-2}{7-2c}$. But $\ff(t)$ is
also equal to $\nu(7) = \frac{5c-25}{3c-21}$. Equating the two leads to the quadratic
equation $13 \, c^2-112 \, c+217=0$, and hence two values of
$c$. Since $f = \nu^{-1}(-1) = \frac{c+10}{8-c}$ is completely
determined by $c$, we get two pairs $(c,f)$. 

The crucial difference between the two examples is that 
three of the elements in $\{A, C, D, F\}$ are specified in the first, and
only two in the second. 
\subsection{} 
We define an assignment to be a bijection $\beta: \FF \ra Z$, for
some $4$-element subset $\FF \subseteq \ltr$. 
An extension of $\beta$ is an alignment $\beta': \ltr \ra
\Gamma$ such that $\beta'|_\FF = \beta$. Here $\Gamma$ is necessarily 
an $R$-configuration containing $Z$. Define the type of $\beta$ to be the cardinality of the set
$\hltr \cap \FF$. 

For instance, the first example corresponds to the assignment 
\begin{equation} 
\bA \ra 2, \quad \bC \ra 3, \quad \bD \ra 5, \quad \bE \ra 7, 
\label{assign} \end{equation} 
which is of type $3$. It admits a unique extension 
\[ \bB \ra \frac{95}{31}, \quad \bF \ra \frac{13}{5}. \] 

\begin{Proposition} \rm 
An assignment has respectively $2, 1$ or $0$ extensions according to
whether its type is $2, 3$ or $4$. 
\end{Proposition} 
\proof Assume that the type is $4$. But since a general quadruple $Z$ is not harmonic, 
it cannot occupy the positions $\{A, C, D, F \}$ in an
$R$-configuration. Hence there cannot be any extensions. 

The first example in section~\ref{examples.type}  illustrates type $3$, and 
the second illustrates type $2$. The proofs in the general case are
exactly on the same lines, hence we leave them to the reader. The only
issue which perhaps requires comment is the following: if the type is
$2$, then we get a quadratic equation for one of the unknown
letters. Since the $z_i$ are general, the equation has two distinct
roots rather than a repeated root. Either of the roots determines the
other unknown uniquely. \qed 

The geometry of the type $3$ case is
utterly straightforward. If three letters from $\hltr$ are specified,
so is the fourth. This fixes Diagram~\ref{fig:ricochet2}, and then
specifying either $B$ or $E$ also specifies the other. 

As to a type $2$ case, assume that $\{A, B, C, E\}$ are
specified. Then so are $V, Z$ and hence the line $ZE$ is specified (on which $W$
must lie). Since $V, D, F$ are collinear, knowing $D$ is tantamount to
knowing $F$. Now, for a variable point $D$ on the conic, the 
function $D \lra AD \cap CF$ traces a conic in the
plane. It intersects $ZE$ in two points, which are the two
acceptable positions of $W$. The other type $2$ cases are similar.

\subsection{} 
Fix a general quadruple $Z$. It takes an elementary counting
argument to see that there are 
\begin{itemize} 
\item $144$ assignments of type $2$, 
\item $192$ assignments of type $3$, and 
\item $24$ assignments of type $4$. 
\end{itemize} 
For instance, to form a type $2$ assignment, choose two letters from
$\hltr$ in $6$ ways. Those, combined with $\{\bB, \bE\}$, can be
distributed in $24$ ways over the $z_i$. Hence there are $24 \times 6 = 144$ 
such assignments. 

Now observe that the group $H$ will act on the sets of
assignments and extensions. For instance, the element 
$v$ in (\ref{defns.uv}) will change the assignment in~(\ref{assign}) to 
\[ \bD \ra 2, \quad \bF \ra 3, \quad \bA \ra 5, \quad \bB \ra 7, \] 
and its extension to 
\[ \bE \ra \frac{95}{31}, \quad \bC \ra \frac{13}{5}. \] 
Of course, both assignments lead to the same $R$-configuration, namely 
$\left\{2, 3, 5, 7, \frac{95}{31}, \frac{13}{5} \right\}$. 

Since elements of $H$ preserve the harmonic subset $\hltr$, they do not affect the
type of an assignment. If two assignments $\beta_1, \beta_2$ are in the same $H$-orbit, then 
the $R$-configurations obtained by extending them will be the
same. Conversely, suppose that $\beta_1, \beta_2$ are two assignments with
respective extensions $\beta_1', \beta_2'$ such that 
$\Gamma = \text{image}(\beta_1') = \text{image}(\beta_2')$. But 
since $Z$ is general, the symmetry group of $\Gamma$ is exactly
$H$, and no larger. Hence $\beta_1, \beta_2$ must be in the same $H$-orbit. 

Now we can count the number of possible $R$-configurations which extend
a given $Z$. There are $\frac{144}{8} \times 2 = 36$ configurations
coming from all assignments of type $2$, and $\frac{192}{8} \times 1 = 24$ 
from those of type $3$. There are none coming from assignments of type
$4$, which gives a total of $24+36=60$. This proves the
following: 

\begin{Theorem} \rm  
The degree of $\Rico$ is $60$. \qed 
\end{Theorem} 

Now we will look for equations which define the variety $\Rico$. The
simplest situation would be that of an ideal-theoretic complete
intersection; i.e., $\Rico$ would be defined by two equations of
degrees $d_1, d_2$ such that $d_1 \, d_2 = 60$. As we will see, this
is not too good to be true. 

\section{Equivariant equations for the ricochet locus} 
\label{section.equiv.equations} 
We begin with a short introduction to classical invariant
theory, which should motivate some of the calculations to follow. The
crucial notion is that of a covariant of binary forms. The
most readable classical references on this subject are~\cite{Elliott,GraceYoung, Salmon}. Modern
accounts may be found in~\cite[Appendix B]{Hunt} and~\cite[Ch.~4]{Sturmfels}.

\subsection{Invariants and Covariants} 
The invariant theory of binary quartics is as good an illustration as
any. Consider a degree $4$ polynomial 
\[ \Phi = a_0 \, x_1^4 + a_1 \, x_1^3 \, x_2 + a_2 \, x_1^2 \, x_2^2 +
a_3 \, x_1 \, x_2^3 + a_4 \, x_2^4,  \qquad (a_i \in \complex)\] 
in the variables $\ux = \{x_1, x_2 \}$. 
Its Hessian, which we denote by $\He(\Phi)$, is defined to be the self-transvectant $(\Phi,
\Phi)_2$. It has an expression 
\[ \He(\Phi) =  \left(\frac{1}{3} \, a_0 \, a_2 - \frac{1}{8} \, a_1^2
  \right) \, x_1^4+ 
\left(a_0 \, a_3- \frac{1}{6} \, a_1 \, a_2\right) \, x_1^3 \, x_2 + \dots + 
\left(\frac{1}{3} \, a_2 \, a_4- \frac{1}{8} \, a_3^2 \right) \, x_2^4. \] 
It is an example of a covariant, since its construction is compatible with a linear change of variables in the
following sense. Given a matrix with determinant $1$, say
$\left[ \begin{array}{cc} 2 & 3 \\ 5 & 8 \end{array} \right]$, we have
the corresponding change of variables: 
\begin{equation} x_1 \lra 2 \, x_1 + 5 \, x_2, \qquad x_2 \lra 3 \, x_1 + 8 \,
x_2. \label{eq.chvar} \end{equation} Now consider the following two
processes: 
\begin{itemize} 
\item Use (\ref{eq.chvar}) in $\Phi$ to get another polynomial
  $\Phi'$, and take its Hessian $\He(\Phi')$. 
\item 
Use (\ref{eq.chvar}) in $\He(\Phi)$ to get $[\He(\Phi)]'$. 
\end{itemize} 
The outcomes are identical, i.e., $\He(\Phi') =
[\He(\Phi)]'$. Since the Hessian is of degree $2$ in the $a_i$, and
degree $4$ in the $\ux$, it is called a covariant of degree-order
$(2,4)$. A covariant of order $0$, i.e., one which contains no
$\ux$-terms, is called an invariant. For instance, 
\[ (\Phi, (\Phi, \Phi)_2)_4 = 
a_0 \, a_2 \, a_4-\frac{3}{8} \, a_1^2 \, a_4-\frac{3}{8} \, a_0 \,
a_3^2+ \frac{1}{8} \, a_1 \, a_2 \, a_3- \frac{1}{36} \, a_2^3, \] 
is an invariant of degree $3$. It is a foundational theorem in the
subject that every covariant is expressible as a compound
transvectant; that is to say, it can be written as a linear
combination of terms of the form 
\[ (\dots (\Phi, (\Phi, \Phi)_{r_1}))_{r_2}, \dots)_{r_k}.  \] 

Any invariant of binary quartics is a polynomial in the two fundamental invariants 
$(\Phi,\Phi)_2$ and $(\Phi, (\Phi,\Phi)_2)_4$. A similar statement is true 
of covariants, but the corresponding list is longer. 

\subsection{} 
The expression $\He(\Phi)$ is identically zero, if and only if $\Phi$ is the fourth
power of a linear form. This illustrates the principle that any property of
a polynomial which is stable under a change
of variables is equivalent to the vanishing of a finite number
of covariants.\footnote{This can be made precise as follows: The space
  of $\P^m$ of binary $m$-ics has coordinate ring $S = \complex[a_0,
  \dots, a_m]$. The action of $SL(2)$ endows $S$ with the structure of a graded
  representation. The locus of polynomials which satisfy a certain invariant
  property is an $SL(2)$-stable subvariety $X \subseteq
  \P^m$, whose ideal $I_X \subseteq S$ is a subrepresentation. Since
  $S$ is a noetherian ring, we can choose a finite number of
  covariants whose coefficients generate this ideal.} 

As another illustration, if $\alpha_i \, x_1 + \beta_i \, x_2,
(i=1,\dots,4)$ are the linear factors of $\Phi$, then 
$(\Phi, (\Phi, \Phi)_2)_4$ is identically zero exactly when the four points 
$[\alpha_i, \beta_i] \in \P^1$ are harmonic, i.e., their cross-ratio
in some order is $-1$. 

\subsection{} 
All of this carries over to polynomials of arbitrary degree $d$, but
the size and complexity of the minimal set of covariants (the so-called `fundamental
system') grow rapidly with $d$. Our immediate interest lies in the case
$d=6$, where the fundamental system has a total of five invariants, namely one each in
degrees $2, 4, 6, 10, 15$.  We will denote them by $I_2, I_4$
etc. Explicit transvectant expressions for the $I_r$ are given
in~\cite[p.~156]{GraceYoung}, but we will not reproduce them here. 

Now, to return to the subject of $R$-configurations, we are looking for covariants which
vanish on the binary sextic 
\begin{equation} 
G_t = \underbrace{x_1 \, x_2 \, (x_1 - x_2) \, (x_1 + x_2)}_\Theta \, 
\underbrace{(x_1 - t \, x_2) \, (x_1 - \ff(t) \, x_2)}_{\Delta_t}. 
\end{equation} 
There is no general procedure which is assured to solve such a  
problem.  However, let us take two plausible decisions at the outset: 
\begin{enumerate} 
\item 
It will be easier to look for invariants, rather than arbitrary covariants. 
\item 
The decomposition $G_t = \Theta \, \Delta_t$ is likely to be helpful,
especially since $(\Theta, (\Theta,\Theta)_2)_4=0$ by the
harmonicity of $\Theta$. 
\end{enumerate} 
These decisions will eventually be vindicated by the fact that they 
lead to a complete solution. Had this not happened, one
would have to start anew and try another strategy. There is no prior guarantee of success. 

\subsection{} Each invariant of $G_t$ is
expressible\footnote{This would technically be true of any sextic
  form, but there is nothing to be gained by chopping up an arbitrary sextic into a
  quartic and a quadratic. This is worth doing here precisely because
  $\Theta$ and $\Delta_t$ are simpler than in the general case.} 
as a polynomial in 
\begin{enumerate} 
\item the individual invariants of $\Theta$ and $\Delta_t$, together
  with 
\item joint invariants of $\Theta$ and $\Delta_t$. 
\end{enumerate} 
The individual invariants are 
\[ \begin{array}{ll} \theta_{20} = (\Theta, \Theta)_4, \; \theta_{30} =
(\Theta, (\Theta,\Theta)_2)_4 \quad \text{for $\Theta$;} \\ 
\delta_{02} = (\Delta_t,\Delta_t)_2, \quad \text{for
     $\Delta_t$.} \end{array} \] 
The joint ones are 
\[ \beta_{12} = (\Theta,\Delta_t^2)_4, \quad 
\beta_{22} = (H, \Delta_t^2)_4, \quad 
\beta_{33} = (T, \Delta_t^3)_6, \] 
where $H = (\Theta,\Theta)_2$ and $T = (\Theta,H)_1$. 
This list is taken from~\cite[p.~168]{GraceYoung}. 
The notation is such that if $\oslash$ stands for any of the letters $\theta, \delta,
\beta$, then $\oslash_{ij}$ is of degree $i$ in $\Theta$ and $j$ in
$\Delta_t$. 

In our case we have 
$\theta_{20} = \frac{1}{2}$, and $\theta_{30}=0$. The remaining invariants are also easy
to calculate; they are as follows: 
\begin{equation} \begin{array}{llll} 
\delta_{02} = - \frac{1}{2} \, \frac{(t^2+1)^2}{(t+1)^2}, & 
\beta_{12} =  \frac{1}{2} \, \frac{(t^2+2t-1) \,
                                                           (t^2-2t-1)}{(t+1)^2},
     & 
\beta_{22} = \frac{\delta_{02}}{3}, & 
\beta_{33} = - \frac{1}{4} \, \frac{t(t-1)(t^2+1)}{(t+1)^2}. 
\end{array} \label{theta.beta.formulae} \end{equation} 
This implies that $\beta_{33}^2 = \frac{1}{32} (\delta_{02}
\beta_{12}^2 - \delta_{02}^3)$, and hence $\beta_{22}, \beta_{33}^2$
are, in effect, redundant. All of this simplifies things considerably. 

\subsection{} The actual derivation of the formulae for $I_r$ needs the symbolic
calculus, as explained in~\cite{Glenn} or~\cite{GraceYoung}. 
Such calculations are tedious and often unpleasant to read through, 
hence we will sketch the derivation of
$I_2$ as an example, and leave the rest as an exercise for the
patient reader. 

We will follow the recipe of \cite[\S 3.2.5]{Glenn}. 
Write $\Theta = a_\ux^4 = b_\ux^4$, and $\Delta_t = p_\ux^2 = q_\ux^2$, where
$a, b, p, q$ are symbolic (or \emph{umbral}) letters. Then 
$I_2 = (a_\ux^4 \, p_\ux^2, b_\ux^4 \, q_\ux^2)_6$ is a sum of $6! = 720$
terms, which are of three kinds: 
\begin{itemize} 
\item $48$ terms of the form $(a \, b)^4 \, (p \, q)^2$,
\item $288$ terms of the form $(a \, b)^2 (a \, q)^2 \, (p \, b)^2$,
  and 
\item $384$ terms of the form $(p \, q) \, (a \, q) \, (p  \, b) (a \,
  b)^3$. 
\end{itemize} 
We have identities 
\[ (a \, b)^4 \, (p \, q)^2 = \theta_{20} \, \delta_{02}, \quad 
(a \, b)^2 (a \, q)^2 \, (p \, b)^2 = \frac{1}{3} \, \theta_{20} \,
\delta_{02} + \beta_{22}, \quad 
(p \, q) \, (a \, q) \, (p \, b) (a \,  b)^3 = \frac{1}{2} \, 
\theta_{20} \, \delta_{02}. \] 
The first is immediate from the definition. The second and the third
follow by a straightforward expansion after using the Pl{\"u}cker  syzygy 
$(a \, q) \, (p \, b) = (a \, p) \, (q \, b) + (a \, b) \, (p
\,q)$. 
And then, 
\[ I_2 = \frac{1}{720} \left[ (48 + 288/3+384/2) \, \theta_{20} \,
\delta_{02} + 288 \, \beta_{22} \right] = 
\frac{7}{15} \, \theta_{20} \, \delta_{02} + \frac{2}{5} \,
\beta_{22}, \] 
which is the required formula. As it stands, it is applicable to any binary sextic written 
as a product of a quartic and a quadratic. But now we can use the 
simplifications in (\ref{theta.beta.formulae}) to get 
\begin{equation} I_2 = \frac{11}{30} \, \delta_{02}. 
\label{I2.formula} \end{equation}

\subsection{} 
With rather more work of the same kind, one deduces the following 
formulae\footnote{The expression for $I_{15}$ is very
  intricate. We can afford to omit it here, because 
  it won't be needed in this calculation.} for the remaining
invariants: 
\[ \begin{array}{ll} I_4 = \frac{2}{1125} \, \beta_{12}^2 + \frac{124}{5625} \, 
\delta_{02}^2, & 
I_6 = \frac{91}{253125} \, \delta_{02} \, \beta_{12}^2 + 
\frac{98}{1265625} \, \delta_{02}^3, \\ 
I_{10} = \frac{416}{284765625} \, \delta_{02} \, \beta_{12}^4 
+ \frac{1141}{284765625} \, \delta_{02}^3 \, \beta_{12}^2 
- \frac{1372}{7119140625} \, \delta_{02}^5. 
\end{array} \] 
Notice that each $I_r$ is expressible as a polynomial in
only two `variables' $x = \delta_{02}$ and $y = \beta_{12}$. 
It follows that $I_2^3, I_2 \, I_4, I_6$ are linear combinations
of the two-element set $\{x^3, x \, y^2 \}$, and hence must
be linearly dependent. The actual dependency relation is easily found
by solving a set of linear equations; it turns out to be 
\begin{equation} 
\underbrace{4032 \, I_2^3-25025 \, I_2 \, I_4+45375 \, I_6}_{\U_6} =0. 
\label{eq.u6} \end{equation} 
The readers may wish to convince themselves that
a parallel argument gives nothing in degrees $2$ or $4$. 

\subsection{} We can use the same line of argument to find another
such invariant in degree $10$. (As before, there is nothing new to be found in degree $8$.)  
The space of degree $10$ invariants for binary sextics is spanned by the six elements 
\begin{equation} I_2^5, \quad I_2^3 \, I_4, \quad I_2 \, I_4^2, \quad I_2^2 \, I_6,
\quad I_4 \, I_6, \quad I_{10}. 
\label{deg10.inv} \end{equation} 
It is contained in the span of the three-element set $\{x \, y^4, x^3
\, y^2, x^5\}$, and hence there must be three linearly independent 
invariants of degree $10$ vanishing on
$\Rico$. Now $\U_6 \, I_2^2 = \U_6 \, I_4 =0$ accounts for
two of these, which leaves room for a new invariant which is not a multiple of $\U_6$. 
Once again, a routine calculation in linear algebra
shows that one can take it to be 
\[ \underbrace{358278336 \, I_2^2 \, I_6-2772533775 \, I_4 \, I_6+
6933745 \, I_2 \, I_4^2+1207483200 \, I_{10}}_{\U_{10}} = 0. \] 
We have arrived at the following statement: 
\begin{Proposition} \rm 
For an arbitrary $t$, we have 
$\U_6(G_t) = \U_{10}(G_t) =0$. \qed 
\end{Proposition}

Let $\Y \subseteq \P^6$ tentatively denote the $4$-dimensional variety defined by the equations $\U_6 =
\U_{10}=0$. By B{\'e}zout's theorem, $\Y$ has 
degree $6 \times 10=60$. Now $\Rico \subseteq \Y$ by the proposition, and since
they have the same degrees, we must have $\Rico = \Y$. We have proved the following: 
\begin{Theorem} \rm 
Let $\Phi$ be a binary sextic representing a set of six distinct
points $\Gamma \subseteq \conic$. Then $\Gamma$ is in $R$-configuration, if and
only if $\U_6(\Phi) = \U_{10}(\Phi) = 0$. \qed 
\end{Theorem} 

We have $I_{\Rico} = (\U_6, \U_{10})$, i.e., $\Rico$ is an ideal-theoretic
complete intersection. This implies that any covariant which vanishes
on $\Rico$ is expressible in the form $f \, \U_6 + f' \, \U_{10}$ for
some $f, f'$. Hence there is no such essentially new
covariant remaining to be found. 

Thus we have completely succeeded in finding invariant-theoretic necessary and
sufficient conditions which characterise the $R$-configuration. A
large part of the success is owed to the fact that the
product of degrees of our two invariants turned out to be exactly
the degree of $\Rico$. In this we have been fortunate, to the extent that such a term has any
meaning in mathematics.

\medskip 

\centerline{--} 

\vspace{1cm}

\parbox{7cm}{ \small 
Jaydeep Chipalkatti \\
Department of Mathematics \\ 
University of Manitoba \\ 
Winnipeg, MB R3T 2N2 \\ 
Canada. \\ \\
{\tt Jaydeep.Chipalkatti@umanitoba.ca}}

\end{document}